\documentclass[12pt]{article}
\usepackage{amssymb,amsmath,amsthm,amscd}

\title{{On a Problem of Bremermann Concerning Runge Domains}
\thanks{ This work was supported by a Marie Curie International Reintegration Grant}}
\author{
Cezar Joi\c ta }
\begin{document}

\newtheorem{df}{Definition}
\newtheorem{lem}{Lemma}
\newtheorem{teo}{Theorem}
\newtheorem{prop} {Proposition}
\newtheorem{cor}{Corollary}
\newcommand{\N}{\mathbb{N}}
\newcommand{\Z}{\mathbb{Z}}
\newcommand{\C}{\mathbb{C}}
\newcommand{\R}{\mathbb{R}}
\newcommand{\nn}{\newline\noindent}
\date{}
\maketitle
\begin{abstract}
In this paper we give an example of a bounded Stein domain in $\C^n$, 
with smooth boundary,
which is not Runge and whose intersection with every complex line
is simply connected.
\end{abstract}
\section{Introduction}

In \cite{Br} Bremermann asked the following question: 

\

"Suppose that
$D$ is a Stein domain in $\C^n$ such that for every complex line
$l$ in $\C^n$, $l\setminus D$ is connected. Is it true that $D$ is
Runge in $\C^n$?" 

\

The question remained open and was mentioned again
in a recent book by T. Ohsawa (\cite{Oh}, page 81). In this paper we will give a negative answer to Bremermann's question. Namely, we will give an example of a bounded, strictly pseudoconvex domain in $\C^n$ with real analytic boundary which
is not Runge in $\C^n$ but whose intersection with every complex line is simply connected.

Note that if $D$ is bounded the hypothesis of the problem means simply that for every complex line $l$, $l\cap D$ is Runge in $l$. If, in addition, one requires that $l\cap D$ is connected as well
then it does follow that $D$ is Runge. See for example \cite{Ho},
page 309, Theorem 4.7.8.

For simplicity our construction will be done in $\C^2$ but it can be 
easily adapted to $\C^n$ for $n\geq 2$. To produce our example we will 
construct first a bounded, strictly pseudoconvex domain $W\subset 
\C^2$ with smooth, real analytic boundary which is Runge but its 
closure is not holomorphically convex. (Note that this is not possible 
in $\C$.)
Next we show that, in fact, we can construct $W$ as above and moreover 
it has
the following geometric property: for every complex line $l$ the set 
of points where $l$
is tangent to $\partial W$ is at most finite. If this is 
the case, then one can show that $l\cap \overline W$ is polynomially 
convex, again for every complex line $l$. Finally, we show that an 
appropriate neighborhood of $\overline W$ is a counterexample to 
Bremermann's problem.

\section{The Example}
The construction will be done in several steps.

First we prove that there exists a bounded domain in $\C^2$ with 
smooth, real analytic
boundary which is strictly pseudoconvex, Runge in $\C^2$, and its 
closure is not polynomially convex.

J. Wermer \cite{We} proved that there exists a biholomorphic map $F$ from a 
polydisc $P=\{(z_1,z_2)\in\C^2:\ |z_1|<a, |z_2|<b\}$ in $\C^2$ onto an open set $F(P)$ of $\C^2$ such that $F(P)$ is 
not polynomially convex. (Wermer's original result was in $\C^3$ but 
it can be
modified to hold in $\C^2$ as well; see \cite{Fo} or \cite{Oh}.) We 
start with such a map and let 
$U_n:=\{z\in\C^2:\ |\frac{z_1}{a}|^n+|\frac{z_2}{b}|^n<1\}$.
Since $U_n\subset U_{n+1}\subset P$ and $\cup U_n=P$, it follows that 
there
exists $m\in\N$ such that $F(U_{m})$ is not polynomially convex. Set
$U:=U_{m}$ and $V=F(U)$. If we define $\varphi:U\to\R$ by $\varphi(z)=
\dfrac 
1{1-|z_1/a|^{2m}-|z_2/b|^{2m}}+|\dfrac{z_1}a|^2+|\dfrac{z_2}b|^2-1$ 
then $\varphi$
is a strictly plurisubharmonic real analytic function and has only one 
critical
point. Since $F$ is a biholomorphism, $\varphi\circ F^{-1}:V\to\R$ has 
the same
properties and it is an exhaustion function for $V$. For $\alpha>0$ 
let $V_\alpha=\{z\in
V:\ \varphi\circ F^{-1}(z)<\alpha\}$. It follows that there exists 
$\alpha>0$ such that
$V_\alpha$ is not polynomially convex. On the other hand, if 
$z_0=F(0)$ (this is the
minimum point and the only critical point of $\varphi\circ F^{-1}$ and
$\varphi\circ F^{-1}(z_0)=0$ ) and we choose $B\subset V$ a ball 
centered at
$z_0$, then there exists $\alpha>0$ such that $V_\alpha\subset B$. 
It follows  that
$V_\alpha$ is Runge in  $B$ (because $\varphi\circ F^{-1}$ is defined 
on $B$) and
therefore is polynomially convex. Put 
$r:=\sup\{\alpha\in\R:\ V_\alpha\ \text{is\
polynomially convex}\}$. From the above observations we deduce that $0<r<\infty$.

We claim that $V_r$ is the example that we are looking for. Indeed $V_r$ is Runge in $\C^2$
as an increasing union of Runge domains and it has smooth, real analytic boundary
because $\varphi\circ F^{-1}$ has no critical point on the $\partial V_r$. We only need 
to convince ourselves that $\overline V_r=\{z\in
V:\ \varphi\circ F^{-1}(z)\leq r\}$ is not polynomially convex.
If $\overline V_r$ was polynomially convex then it would have a Runge (in $\C^2$)
neighborhood $W$ with $W\subset V$. If this was the case then for $\epsilon> 0$
small enough $V_{r+\epsilon}\subset W$ and $V_{r+\epsilon}$ would be Runge
in $W$ and therefore in $\C^2$. This would contradict the choice of $r$.

Let us refrase what we have done so far. We proved that if $V$ is a domain in $\C^2$ and $\phi:V\to \R$ is
a strictly plurisubharmonic function such that there exist $a_0<a_1$
real numbers with the following properties:
$$\left.\begin{array}{l}
\{x\in V:\phi(x)<a_1\}\subset\subset V, \cr
\{x\in V:\phi(x)<a_0\}\ {\text {is connected and
 contains}}\ C(\phi):=
{\text{the set of }}\cr
{\text{critical points of }} \phi,\cr
\{x\in V:\phi(x)<a_0\}\ {\text{is Runge in }}\C^2,\cr
\{x\in V:\phi(x)<a_1\}
\ {\text{is not Runge in }}\C^2
\end{array}
\right\}(*)$$

Then there exists a unique real number $r=r(\phi)\in[a_0,a_1)$ such that 
$V_{r(\phi)}:=\{x\in V:\phi(x)<r(\phi)\}$ is Runge and $\overline V_{r(\phi)}=
\{x\in V:\phi(x)\leq r(\phi)\}$ is not holomorphically convex. Note that
$V_{r(\phi)}$ must be connected since each of its components
contains a (minimum) critical point, $V_{r(\phi)}$ contains $\{x\in V:\phi(x)<a_0\}$ which is connected and $\{x\in V:\phi(x)<a_0\}\supset C(\phi)$. We also proved that there exists a real analytic function
$\phi$ satisfying (*). We fix such a $\phi$. Shrinking $V$ we can assume that $\overline V$ is compact and that $\phi$ is defined
on a neighborhood of $\overline V$.  

\

Next we want to show that there exists $\psi$, a small perturbation of $\phi$,
which satisfies (*) and in addition it has the following geometric property:
for every complex line $l$ the set ${\cal T}(\psi,l):=\{x\in V_{r(\psi)}\cap l:l\ {\text{is tangent to}}\ \partial V_{r(\psi)}
\ {\text {at}}\ x\}$ is finite. 

Indeed: let $U$ be an open and connected set such that $C(\phi)\subset
U\subset\subset V_{r(\phi)}$ and let $W$ be an open and
relatively compact neighborhood of $\partial V_{r(\phi)}$ 
and $0<\delta<\delta'<a_1-r(\phi)$ two real numbers such that $U\subset
\subset \{x\in V:\phi(x)<r(\phi)-\delta\}$ and
$\{x\in V:r(\phi)-\delta<\phi(x)<r(\phi)+\delta\}\subset\subset W
\subset\subset \{x\in V:\phi(x)<r(\phi)+\delta'\}$.

If $\epsilon>0$ is small enough then, for every $\psi:V\to\R$, a ${\cal C}^\infty$ function,
if the $\sup$ norms on $\overline V$ of $\psi-\phi$, $\frac{\partial(\psi-\phi)}
{\partial x_j}$, $\frac{\partial^2(\psi-\phi)}{\partial x_j\partial x_k}$,
$j,k=1,...,n$ (here we denote $z_j=x_{2j+1}+ix_{2j}$) are less than $\epsilon$
then $\psi$ is strictly plurisubharmonic and satisfies (*). Moreover 
$C(\psi)\subset U$, $r(\psi)\in[r(\phi)-\delta,r(\phi)+\delta]$ and
$\partial\{x\in V:\psi<s\}\subset W$ for every 
$s\in[r(\phi)-\delta,r(\phi)+\delta]$. We claim that there exists a real analytic
$\psi$ such that for every complex line $l$ and for every 
$s\in[r(\phi)-\delta,r(\phi)+\delta]$ the following set 
${\cal T}(\psi,l,s):=\{x\in\partial
\{x\in V:\psi(x)<s\}\cap l:l\ {\text{is tangent to}}\ \partial\{x\in V:\psi(x)<s\}
\ {\text {at}}\ x\}$ is finite. Indeed, if $x_0$ is not isolated in ${\cal T}(\psi,l,s)$
and if we denote by $u:=\psi_{|l}$ then $x_0$ is not isolated in 
$\{z\in l: u(z)=u(x_0),\ \nabla u(z)=0\}$. On the other hand
$u$ is real analytic and strictly subharmonic. It follows that around
$x_0$ at least one of the sets $\{z\in l\cap V:\frac{\partial u}{\partial x}(z)
=0\}$ or $\{z\in l\cap V:\frac{\partial u}{\partial y}(z)
=0\}$ is smooth and then the smooth one is contained in the other
one. Hence there exists around $x_0$ a smooth real analytic curve $C$ such that $u_{|C}=u(x_0)$ and $\nabla u_{|C}=0$. If $\{f=0\}$
is a local equation for $C$ it follows that $\psi_{|C}=u=u(x_0)+
f^2g$. However it is not difficult to see that this condition is
not satisfied by a generic real analytic function $\psi$.
(For example one notices that $\det Hess (u)(x_0)=0$
and after a linear change of coordinates we can assume
that $f(x)=x+$ higher order terms and $g(x_0)=1$. Then a 
straightforward computation shows that $u$ must satisfy the following
four conditions at $x_0$: 
$\dfrac{\partial^3u}{\partial y^3}=0$,
 $\dfrac{1}{4!}\dfrac{\partial^4u}
{\partial y^4}=[\dfrac{1}{4}\dfrac{\partial^3 u}{\partial x\partial y^2}]^{2}$,
$\dfrac{1}{5!}\dfrac{\partial^5u}{\partial y^5}=[
\dfrac{1}{3!}\dfrac{\partial^4 u}{\partial x\partial y^3}-
\dfrac{1}{8}\dfrac{\partial^3 u}{\partial x\partial y^2}\dfrac{\partial^3 u}{\partial x^2\partial y}]\dfrac{1}{4}
\dfrac{\partial^3 u}{\partial x\partial y^2}$. \Big)

\ \nn
We fix now a $\psi$ which satisfies (*) and the geometrical
property mentioned above. 

\

Our next goal will be to show that, 
for every complex line $l$ in $\C^2$, $l\cap \overline V_{r(\psi)}$ is polynomially convex (although $\overline V_{r(\psi)}$ is not). 
Note that $(l\cap\overline V_{r(\psi)})\setminus \overline{l\cap V_{r(\psi)}}$ is 
a finite set (as a subset of ${\cal T}(\psi,l)$).
Hence it suffices to show that $\overline {l\cap V_{r(\psi)}}$ is polynomially convex. Let's assume that it is not. Note
that $l\cap V_{r(\psi)}$ is Runge in $l$ (since $V_{r(\psi)}$ is Runge in $\C^2$)
and that it has a smooth boundary except at a finite set of points (the set of points of non-smoothness is also a subset of ${\cal T}(\psi,l)$). As we assumed that $\overline{l\cap V_{r(\psi)}}$ is not polynomially convex it follows that there exists
a rectifiable loop $\gamma$ in $l$ such that $\gamma\setminus(l\cap V_{r(\psi)})$
contains only points where the boundary of $l\cap V_{r(\psi)}$ in $l$
is not smooth and therefore is finite and $\widehat\gamma\cap(l\setminus(\overline{l\cap V_{r(\psi)}}))\neq\emptyset$ (in fact it has a nonempty interior). Using again the finiteness of ${\cal T}(\psi,l)$ it follows that
$\widehat\gamma\cap(\C^2\setminus\overline V_{r(\psi)})\neq\emptyset$.
We claim that there exists a $\cal C^\infty$ family of biholomorphisms $\{f_\epsilon:\C^2\to\C^2\}_{\epsilon\in\R}$
such that $f_0$ is the identity and for $\epsilon>0$ small enough
$f_\epsilon(\gamma)\subset V_{r(\psi)}$. Without loss of generality
we can assume that $l=\{z=(z_1,z_2)\in\C^2:z_2=0\}$.
We write $\gamma\setminus(l\cap V_{r(\psi)})=:\{(p_1,0),\dots,(p_s,0)\}$
and we denote by $(0,q_1),\dots(0,q_s)$ the unit inner normals to
$\partial V_{r(\psi)}$. We choose $h:\C\to\C$ a holomorphic function
such that $h(p_j)=q_j$ and we define $f_\epsilon(z)=(z_1,z_2+\epsilon h(z_1))$. It is obvious that $f_\epsilon$ are biholomorphisms and since $\dfrac{{\rm d}f_\epsilon}{{\rm d}\epsilon}(p_j,0)=(0,q_j)$ it follows that $f_\epsilon$ have the sought properties. Because 
$\widehat\gamma\cap(\C^2\setminus\overline V_{r(\psi)})\neq\emptyset$ and $\{f_\epsilon\}$ is a continuous family we deduce  that
for $\epsilon$ small enough $f_\epsilon(\widehat\gamma)\not\subset \overline V_{r(\psi)}$. On the other hand $f_\epsilon(\widehat\gamma)=\widehat{f_\epsilon(\gamma)}$
and $f_\epsilon(\gamma)\subset V_{r(\psi)}\cap f_\epsilon(l)$. It 
follows from here that $V_{r(\psi)}\cap f_\epsilon(l)$ is not Runge in 
$f_\epsilon(l)$ which is a contradiction since $V_{r(\psi)}$ is Runge 
in $\C^2$ and $f_\epsilon(l)$ is a closed analytic submanifold
in $\C^2$. 

\

We are now ready to produce our example. For $\epsilon>0$
we set $W_\epsilon:=\{x\in V: \psi(x)<r(\psi)+\epsilon\}$. 
It follows from the definition of $r(\psi)$ that $W_\epsilon$
is not Runge in $\C^2$. 
We wish to prove that there exists $\epsilon>0$ such that
for every complex line $l$, $W_\epsilon\cap l$ is Runge in $l$. 

Suppose that this is not the case. Then for $n\in\N$ large enough there exists a complex line $l_n$ such that $W_{\frac 1n}\cap l_n$ is not Runge in $l_n$.Note that $\{l_n\}$ is a sequence of lines that intersect a given compact subset of $\C^2$. It contains then a convergent subsequence. By passing to this subsequence we can assume that $\{l_n\}$ converges to a line $l$. 

We already proved that $l\cap
\overline V_{r(\psi)}$ is holomorphically convex and this implies that there exists
$\Omega$ a Runge open subset of $\C^2$ such that $l\cap
\overline V_{r(\psi)}\subset\Omega\subset V$. As $\cap W_{\frac 1n}=\overline V_{r(\psi)}$ and $l_n$ converges to $l$ we deduce that there exists $n_0\in\N$ such that for every $n\geq n_0$
$W_{\frac 1n}\cap l_n\subset \Omega$. Hence $W_{\frac 1n}\cap l_n=
(W_{\frac 1n}\cap\Omega)\cap l_n$. On the other hand, $\psi$ is a
plurisubharmonic function defined on the whole $\Omega$ which is Stein
and therefore $W_{\frac 1n}\cap\Omega$, which is a level set for $\psi_{|\Omega}$, is Runge in $\Omega$. Since $\Omega$ is Runge in $\C^2$ it follows that $W_{\frac 1n}\cap\Omega$ is also Runge in
$\C^2$ and from here we obtain that $W_{\frac 1n}\cap l_n$ is Runge in
$l_n$. This contradicts our assumption. 

In conclusion, we proved that for $\epsilon>0$ small enough
$W_\epsilon$ is bounded, strictly pseudoconvex, is not Runge in $\C^2$ and for every complex line $l$ in $\C^2$, $W_\epsilon\cap l$ is Runge in $l$. In the same way as before
$W_\epsilon $ must be connected since each of its components contains a critical point of $\psi$.

\ \nn
$\mathbf{Acknowledgments}:$ \* \* {\it  I am very grateful
to Professor Mihnea Col\c toiu for bringing Bremermann's problem to my 
attention and to Professor Terrence Napier for 
very useful discussions.
}

\vspace{1.0cm}
\begin{flushleft}
Cezar Joi\c ta \newline
Institute of Mathematics of the Romanian Academy\newline
P.O. Box 1-764, Bucharest 014700\newline 
ROMANIA\newline
\emph{E-mail address}: Cezar.Joita@imar.ro
\end{flushleft}


\begin{thebibliography}{99}

\bibitem{Br} Bremermann, H.J.: Die Charakterisierung Rungescher Gebiete durch plurisubharmonische Funktionen. {\it Math. Ann.} {\bf 136}, 173--186, 1958. 
\bibitem{Ho} H\"ormander, L.: {\it Notions of convexity.} Progress in Mathematics, 127. Birkh\"auser, 1994. 
\bibitem{Fo}  Forn\ae ss, J.E.; Stens\o nes, B.: {\it Lectures on Counterexamples in Several Complex Variables}. Mathematical Notes, 33. Princeton University Press, 1987.
\bibitem{Oh} Ohsawa, T.: {\it Analysis of several complex variables.}  Translations of Mathematical Monographs, 211.  American Mathematical Society, 2002.
\bibitem{We} Wermer, J.:  An example concerning polynomial convexity.  {\it Math. Ann.} {\bf 139}, 147--150, 1959. 

\end{thebibliography}
 \end{document}